\magnification 1100

\def\irr{1}
\def\split{2}
\def\mult{3}
\def\cusp{4}
\def\splitthm{5}
\def\genus4{6}

\font\tenmsb=msbm10
\font\sevenmsb=msbm10 at 7pt
\font\fivemsb=msbm10 at 5pt
\newfam\msbfam
\textfont\msbfam=\tenmsb
\scriptfont\msbfam=\sevenmsb
\scriptscriptfont\msbfam=\fivemsb
\def\Bbb#1{{\fam\msbfam\relax#1}}

\def\sqr#1#2{{\vcenter{\vbox{\hrule height.#2pt \hbox{\vrule width.#2pt
height #1pt \kern #1pt \vrule width.#2pt}\hrule height.#2pt}}}}
\def\qed{$\  \sqr74$}

 \def\la{\longrightarrow}

 \def\ni{\noindent}
 \def\cl{\centerline}
\def\df{\noindent {\bf Definition.\  }}
 \def\rk{\noindent {\it Remark.\  }}
 
 \def\pf{\noindent {\it Proof.\  }}

\def\d{\delta }

\def\a{\alpha}

\def\g{\gamma}
\def\t{\theta}
\def\T{\Theta}

\font\gothicf=eufm10
\font\sgothic=eufm7
\font\ssgothic=eufm5
\textfont5=\gothicf
\scriptfont5=\sgothic
\scriptscriptfont5=\ssgothic

\def\P{{\Bbb P}}
\def\Pl{\P ^2}

\def\Pic{{\rm Pic}}

\def\Im{{\rm Im}}

\def\Zn{Z^{\nu}}

\cl{\bf{ON MODULAR PROPERTIES OF ODD THETA-CHARACTERISTICS}}

\

\cl{Lucia Caporaso}

\

\

\ni{\bf{Abstract.}} A general canonical curve $X\subset \P ^{g-1}$ determines a
finite set $\t(X)$ of hyperplanes that are tangent to $X$ at $g-1$ points. 
Such a set is in bijective correspondence with the set of odd theta-characteristics
of
$X$. 
We generalize the definition  of $\t(X)$  to certain 
singular curves (including Deligne-Mumford stable curves)
in a way that is compatible with degenerations. 
We give the explicit description of $\t(X)$, with its enumerative and
local data.  This is applied
to show that  some singular curves can be recovered from $\t(X)$.
This paper generalizes the preliminary steps used in [CS] to show that, if $g=3$,
$\t (X)$ uniquely determines $X$.

\

\

\ni{\bf{1. Introduction.}}
A theta-characteristic on a smooth, projective, complex curve $X$ of genus $g$ is a square root of the
canonical line bundle $K_X$, that is, a divisor class $D\in \Pic X$ (of degree $g-1$)
such that $2D = K_X$. Such a $D$
is called ``odd" (respectively ``even") if $h^0(X,D)$ is odd (respectively even).

Assume $X$ to be non-hyperelliptic, fix a canonical model $X^c$  for $X$ and let $H$ be a hyperplane
that is tangent to $X^c$ at $g-1$ (not necessarily distinct) points, that is,
$H\cap X^c = \sum 2p_i$ where $p_i \in X^c$; clearly the divisor $H\cap X^c$ is an effective
theta-characterstic. Conversely, if $D$ is an effective theta-characteristic, then there exists a
hyperplane $H\subset \P ^{g-1}$ such that $D={1\over 2} H\cap X^c$, in particular, $H$ is tangent to
$X^c$ at $g-1$ points. We shall call such hyperplanes the {\it theta-hyperplanes} of $X^c$.

It is well known that a smooth curve has exactly $2^{g-1}(2^g -1)$ odd theta-characteristics.
A simple consequence of the Petris's conjecture
(proved by Gieseker, see [ACGH] V.1.7) is that,
if the
curve is general in moduli, then its odd
theta-characteristics coincide with the effective ones.
In other words a general smooth curve $X$ of genus $g$ has exactly
$$N_g:=2^{g-1}(2^g -1)$$
 effective divisors
$D$ such that
$$2D = K_X\  \   {\rm{and}}\  \  \  \  h^0(X,D) =1.$$
Hence the canonical model of a general $X$ has $N_g$ distinct theta-hyperplanes.

Another intrinsic interpretation can be given using the Jacobian variety $J(X)$
polarized by a  theta divisor $\T _X$: effective theta-characteristics are in natural one-to-one
correspondence with the $2$-torsion points of $J(X)$ contained in $\T _X$. If $X$ is general, then the
Riemann Singularity Theorem implies that $\T _X$ is smooth at such $N_g$ points. As a consequence
the images of the odd theta-characteristics under the Gauss map of $\T _X$ are, once again,
put in bijective correspondence with the theta-hyperplanes of a canonical model of $X$.

Summarizing, the theta-characteristics of a curve $X$ can be viewed abstractly
as the $2$-torsion points of $J(X)$ lying on  $\T _X$. Or projectively, as the hyperplanes that are
tangent to a canonical model of $X$ in $g-1$ points. This second point of view allows one to compare the
odd theta-characteristics of different curves. This paper is written in such a perspective.

Let $X$ and $Y$ be two abstract curves as above. 
If we are given for them canonical models $X^c$ and $Y^c$ in the same $\P ^{g-1}$,
 we can ask: if $X^c$ and $Y^c$ have the same $N_g$ theta-hyperplanes does it follow that $X=Y$?

We shall denote by $
\t (Z)$ the union of all  theta-hyperplanes of a smooth canonical curve $Z$.

Another version of the above question is: given two
distinct canonical curves $Z$ and $Z'$, is it always true that $\t(Z) \neq \t (Z')$?

In [CS]  the question above is answered  affirmatively for curves of genus $3$: 
it is shown that a general plane quartic can be  recovered
by its $28$ bitangent lines.
We conjectured that the same holds for any genus.

The  method in [CS] is a   degeneration technique, 
which can be summarized in three steps as follows.

The first step is to describe how the set of 28 bitangent lines degenerates when smooth curves degenerate
to (certain) singular ones. This requires constructing a geometrically meaningful analogue for $\t (X)$
which  is compatible with deformations to smooth curves, 
that is to say, that does not depend on the choice of the family.

The second step is to show that some singular curves 
can in fact be recovered from their generalized set of bitangent lines
 (constructed in the previous step).

Finally one needs an argument to reduce the proof of the main theorem for smooth curves to the second step.
This is done in [CS] using Geometric Invariant Theory and deformation theory for plane curves.

In this paper we generalize the first two stages of the process to curves of any genus.

For the first step, we concentrate on those singular curves that have been so far proved 
to be the ``best"degenerations of
smooth ones: stable curves in the sense of Deligne-Mumford.
Denote by $\overline {M_g}$ the moduli space of stable curves of genus $g$.

Maurizio Cornalba exhibited in [C] a 
geometrically meaningful compactification of the moduli space of pairs

\cl{(smooth curve $X$, theta-characteristic on $X$).}

\

\ni
His space (called $\overline {S_g}$) is finite over $\overline {M_g}$ 
and has good functorial properties.
He furthermore shows  that $\overline {S_g}$ has two disjoint irreducible components, 
$\overline {S_g}^+$ and $\overline {S_g}^-$,
corresponding to the compactification of the loci of (respectively)
even and odd theta-characteristics of smooth curves.
We shall here give a closer look at $\overline {S_g}^-$ (whose degree over $\overline {M_g}$
 is $N_g$)  describing the enumerative and local data for its natural strata
over the loci of $\overline {M_g}$ parametrizing curves with a given number of nodes.

The second step will be carried out combining the previous analysis with the results of [CS]. 
We prove (in Theorem \splitthm ) that a general, nodal
canonical curve $X\subset \P ^{g-1}$
which is the union of two rational normal curves meeting at $g+1$ points
is uniquely determined by its generalized $\t (X)$ ($X$ is a so called ``split" curve, see below).

Finally we show that singular, irreducible, general  curves of (arithmetic) genus $4$
can be recovered by their odd
theta-characteristics in the following sense: if $X$ and $Y$ are singular, irreducible curves in
$\overline{M_4}$ such that there exist suitable canonical models $X^c$ and $Y^c$ having the same
theta-hyperplanes, then $X=Y$ (this is Theorem \genus4 ).

\

\

\ni
{\bf 2. Enumerative results.} Fix $g\geq 3$, let $r = g-1$ and $d=2g-2$. Denote by
$H_{d,g}^r$ the  Hilbert scheme of curves in
$\P ^r$ having Hilbert polynomial $p(x)=dx-g+1$.
 Throughout $X$ will be a curve in $\P ^r$ satisfying the following set of
assumptions:

\

$X$ is reduced of degree $d$ and arithmetic genus $g$,
$X$ has at most nodes as singularities,
no irreducible
component of
$X$ is contained in a hyperplane,
$X$ is embedded in $\P ^r$ by the (complete) linear series $|\omega _X|$
(where $\omega _X$ denotes the dualizing line bundle of $X$).

\

A curve satisfying the above requirements will be called throughout a {\it canonical} curve.
We shall denote by $V'$ the open subset of $H_{d,g}^r$ parametrizing curves satisfying
the above properties.

Let $\d$ be the number of nodes of $X$; then either $X$ is irreducible and $0\leq \d \leq g$,
or $X$ is reducible, in which case it is easy to see that $X$ must be the union
of 2 rational normal curves (of degree $g-1$ of course) meeting transversally at $g+1$ points.

\

\df
Let $H\subset \P ^r$ be a hyperplane; we shall say that $H$ is a {\it theta-hyperplane} of $X$
if $H\cap X$ is everywhere non-reduced and if $H$ meets $X$ in $g-1$ distinct points.
Let $i$ be an integer with
$0\leq i$. We shall say that $H$ is a theta-hyperplane of type $i$ for $X$ if $H$
contains exactly
$i$ singular points of
$X$.
We shall denote by $\t _i(X)$ the set of such hyperplanes and by $t_i(X)$ its cardinality.

\

As already noticed in the introduction, a  general, smooth $X$  has a
finite number of  distinct theta-hyperplanes, and this number is equal to
$N_g = 2^{g-1}(2^g-1)$; these are the the hyperplanes that
  are tangent to $X$ at $g-1$ points.

We shall say that a smooth
curve $X$ as above is {\it theta-generic} if it has $N_g$
theta-hyperplanes. 

We can define a rational map
$$
\t : V' \cdot \cdot \cdot \cdot \cdot \cdot \ \cdot >  Sym ^{N_g}(\P ^r)^*
$$
by the rule $\t (X) = (H_1,....,H_{N_g})$ with $X$ smooth, theta-generic,
$H_i$ a theta hyperplane of $X$ and $H_i\neq H_j$.

It is not hard to see (arguing like in  [CS] 2.3.1 or using [C]) that $\t$ can be defined 
on an open subset of the Hilbert scheme larger than $V'$, including some singular curves.
More precisely,
if  $X\subset \P ^r$ is a reduced, nodal, canonical curve with $\d$ nodes,
having no component contained in a hyperplane,
and general among all curves with $\d$ nodes,
$X$ has a finite set of distinct theta-hyperplanes. Thus one can define $\t (X)$
which will be a hypersurface of degree $N_g$ all of whose irreducible components are
theta-hyperplanes of $X$, which will now appear with  a  certain multiplicity in 
$\t(X)$.

We shall call such a curve  {\it theta-generic} and we shall denote by $V$ the set of all such curves. Of course
$\t$ can be 
defined over $V$.

To describe the hypersurface $\t (X)$ we start by computing $t_i(X)$.

The multiplicity with
which each hyperplane occurs in $\t(X)$ will be computed later.

\

\proclaim Proposition \irr. Let $X\in V$ be irreducible with $\d$ nodes. Then

$
t_0(X) = 2^{2(g-\d )}2^{\d -1}=2^{2g-\d -1}
$ \  {\it if} \  $\d\geq 1$;

\

$
t_i(X) = 2^{2(g-\d )}2^{\d -i -1} {\d \choose i}={\d \choose i} 2^{2g-\d -i-1}
$ \  {\it if} \  $i< \d< g-1$;

\

$
t_{\d }(X)= 2^{g-\d -1}(2^{g-\d} -1)
$ \  {\it if} \ $\d <g-1$;
$$
t_{g-1}=\cases{1, & if $\d = g-1$;\cr
g, & if $\d = g$.\cr}\hskip3.3in 
$$

\

\rk
We stated the result in such a redundant fashion, with  double equalities, for the reader's
convenience: the first formula  highlights the inductive pattern, while the second is more compact.

\

\pf
We start by proving that, if $\d \geq 1$, then there are $2^{2(g-\d )}2^{\d -1}$ theta-hyperplanes of
type $0$. Let $H$ be such a hyperplane, then 
$$
H\cap X = 2 \sum _{i=1}^{g-1} P_i 
$$
where $P_i$ is a smooth point of $X$. Let $\nu : X^{\nu} \la X$ be the normalization of $X$
and let
$$
E=\sum _{i=1}^{g-1}\nu ^{-1}( P_i).
$$
Now denote by $D$ the class of $\nu ^* H$ in the Picard group of $X^{\nu }$,
so that $D = K_{X^{\nu}} + G$ where $G$ is an effective divisor of degree $2 \d$
and $\deg D = 2g-2$.
Of course $E$ is an effective divisor such that $2E \sim D$; conversely
any divisor class $E'$, such that $2E' \sim D$, is effective, since the genus of $X^{\nu}$ is
$g-\d$ wich is at most equal to $\deg E'$.

Let $W\subset \P(H^0(X^{\nu}, D)^*)$ be the $r$-dimensional linear subspace such that our
ambient space $\P ^r$ is equal to $W$.
Given a square root $E$ of $D$, the number of theta-hyperplanes of type $0$ of $X$ coming
from the linear system $|E|$ is thus equal to the number of divisors in $|E|$ whose square is
in the sublinear system $W$.
Call $x(\d )$ such a number  (we will see that it does not depend on the choice of
$E$). Since on $X^{\nu}$ there are exactly $2 ^{2(g-\d )}$ distinct divisor classes
whose square is equal to
$D$,
we get that
$$
t_0(X)= 2 ^{2(g-\d )}x(\d ).
$$
We now compute $x(\d )$.
There is a natural map
$$
\alpha :|E|\la |D|
$$
sending $E$ to $2E$.
Now, by Riemann-Roch, $\dim |D|= g+\d -2$ and $\dim |E|= \d -1$ (our assumption that $X$ is
theta-generic implies that the divisor $E$ on $X^{\nu}$ is non-special).
Notice  that if $\d = 1$ our result follows immediately.

The map $\a$ is  the composition of a  Veronese mapping of degree $2$
from $\P ^{\d -1} = |E|$, with a general linear projection; thus the degree of
$\Im \alpha$ is equal to $2^{\d -1}$. Since $W$ is a linear subspace of $|D|$ of codimension
$\d -1$ (transverse to $\Im \alpha$ since $X$ is theta-generic),
we conclude that
$$
x(\d )= \deg (W \cap \Im \alpha ) = 2 ^{\d -1}
$$
which is what we wanted.

For the second formula we start with the case $i=1$ (and $\d \geq 2$) for simplicity.
If $i=1$ let $N$ be a fixed node of $X$ and let us compute the number of theta hyperplanes of
type $1$  containing $N$. Project $X$ from $N$ to $\P ^{r-1}$, the image is a nodal
irreducible curve $Y$ having $\d -1$ nodes;  $Y$ is also canonical, not contained
in a hyperplane and theta-generic. The hyperplanes of type 0 of $Y$ are in bijective
correspondence with the hyperplanes of type 1 of $X$  containing $N$.
By the previous  formula, we obtain
$$
t _0(Y) = 2^{2(g-1-\d+1)}2^{\d -2}=2^{2(g-\d)}2^{\d -2}
$$
and, since there are $\d$ nodes, $t_1 (X) = \d \cdot t_0(Y)$;  we are done in this
case.

Let now $1<i<g-1$; we proceed in a similar fashion. Pick a set $S$ of $i$ nodes and
consider first the case $i<g-2$.
Project
$X$ from the $i-1$-dimensional linear space spanned by $S$
(if $S$ spans a set of smaller dimension, pick a general space of dimension $i-1$ containing it;
this is not going to happen for the general $X$); the
image is a canonical curve
$Y\subset \P ^{g-1-i}$ of genus $g-i$ and having $\d -i$ nodes.
The  theta-hyperplanes of type $i$ containing $S$ are in bijective correspondence with the set
of theta-hyperplanes of type $0$ of $Y$. Since there are ${\d \choose i}$ choices
for $S$ we get, using the first formula
$$
t_i(X) = {\d \choose i}t_0(Y)={\d \choose i}2^{2(g-i-\d+i)}2^{\d - i -1}
$$
and we are done.
Now assume $i=g-2$. Let $\pi :X \la X'\subset \Pl$ be the projection from $g-3$ nodes
in $S$ and let $N' \in X'$ be the image of the remaining node in $S$ (the one not
contained in the center of the projection). The curve $X'$ is an integral nodal curve
of degree $4$ having $\d - g +3$ nodes. The set of theta-hyperplanes of type $g-2$ of
$X$ containing $S$ is in bijective correspondence with the set 
$\t _{N'}(X')$ of theta-lines 
of type $1$ of $X'$ containing $N'$.
By induction
$$
\# \t _{N'}(X') ={1\over \d - g +3 }t_1(X') = 2 ^{g-\d}2^{\d - g +1}.
$$
Since there ${\d \choose g-2}$ choices for $S$ we get that
$$
t_{g-2}(X)={\d \choose g-2}\#\t _{N'}(X')={\d \choose g-2}2 ^{g-\d}2^{\d - (g-2) -1}
$$
which is what we wanted.

The formula for $t_{\d}(X)$ is obtained again projecting from the nodes of $X$ and
getting a smooth curve $Y$. The number of theta-hyperplanes of type $\d$ of $X$
is equal to the number of effective theta characteristics of the theta-generic, smooth, canonical curve
$Y$.

The last two equalities
are clear: just notice that $g-1$ nodes on a curve of degree $2g-2$ must span a hyperplane
(or the curve is contained in a hyperplane).
\qed

\

\df
A curve $X$ in $V$ is called a {\it split} curve if $X=C_1\cup C_2$ is the union of two rational normal
curves meeting transversally at $g+1$ points.

\

\rk Split curves are the only reducible curves in $V$, and they are characterized by the property of
having $g+1$ nodes.

\

\proclaim Proposition \split.
Let $X$ be a split curve of (arithmetic) genus $g$. Then for $0\leq j \leq  g-1$ we have:
$$
t_j(X) = \cases{{g+1 \choose j}2^{g-j-1} & if \  $j \not\equiv g $  (mod 2) \cr
0 & if \  $j \equiv g $  (mod 2) \cr}
$$

\

\rk In particular, if $g$ is odd, then $t_1(X)=t_3(X)=....=t_{g-2}(X)=0$. 

If $g$ is even then $t_0(X)=t_2(X)=....=t_{g-2}(X)=0$. 

Notice also that for every $g$, a split curve has $g+1 \choose 2$ theta-hyperplanes of type $g-1$.

\

\pf
Start with $t_0(X)$. If $g$ is even, that is, $r$ is odd,
then 
$t_0(X)$ must vanish: let $X=C_1\cup C_2$ and let $H$ be any theta-hyperplane of type 0.
Then $H$ is tangent to $X$ at $r$ distinct smooth points, since $r$ is odd then $H$ must be tangent to
(say)
$C_1$ in at least $r+1 \over 2$ points, so that $\deg H\cap C_1 \geq r+2$ which is obviously impossible.

Now  consider $g$ odd and write $r = 2h$.
Let $C$ be a rational normal curve in $\P^r$ and let $T_C\subset (\P ^r )^*$ be the closure of the
locus of
hyperplanes in $\P ^r$ that are tangent to $C$ at $h$  distinct points. Then we claim
that $T_C$ is an
irreducible subvariety in $(\P ^r )^*$ of dimension $h$ and degree $2^h$.
The first part of the statement is immediate: just notice that through any $h$ points 
$P_1, ....,P_h$ of $C$ there exists
a unique hyperplane $H_{P_1, ....,P_h}$ that is tangent to $C$ at those $h$ points.
In fact, let $L_i$ be the line tangent to $C$ at $P_i$; then $L_1,....,L_h$ must be in general position, 
since $\deg C = 2h$, hence they span a
unique hyperplane. We get that $T_C$ can be constructed as the closure of the image of the map
$$
\matrix{Sym ^h(C) \setminus {\rm diagonals} &\la &(\P ^r )^*\cr
(P_1,...., P_h)&\mapsto &H_{P_1, ....,P_h}\cr}
$$

Now to compute the degree, let $J\cong \P ^h$ be a general $h$-dimensional
linear subspace of $\P ^r$. The degree   of $T_C$ is
equal to the number of hyperplanes containing $J$ and tangent to $C$ at $h$ points.
Let $\pi :C\la \P ^h$ be the projection from $J$. Having chosen $J$ general, the image $\pi (C) $ is a nodal,
rational curve $D$ of degree $2h$ in $\P ^h$ having $h+1$ nodes.
The theta-hyperplanes of type 0 of $D$ are in bijective correspondence with the hyperplanes 
in $\P ^r$ containing $J$ and tangent to $C$ at $h$ points.
By Proposition \irr , $\  t_0(D) = 2 ^h=2 ^{g-1}$; we conclude that $\deg T_C = 2^h$.

A consequence of this analysis is that if $X=C_1\cup C_2$ is a split curve of odd genus,  the number of
its theta-hyperplanes of type 0 is obtained by intersecting $T_{C_1}$ with $ T_{C_2}$: this is the
intersection of two subvarieties of dimension $h$ and degree $2^h$ in a projective space of
dimension $r=2h$. Transversality is ensured by the condition that $X$ is theta-generic, so that we get
$$
\deg T_{C_1}\cap T_{C_2}= 2^{2h}=2^{g-1}
$$
and we are done.
The rest of the proof of the proposition is by iterating projections of $X$ onto a lower degree split
curve.
Let us show the first case, the others are obtained in the same fashion.

To see that a split curve $X$ of odd genus $g$ has no hyperplanes of type 1, project from a node
$N$ onto a split curve $Y$ of (even) genus $g-1$. Theta hyperplanes of type 1 of $X$ containing $N$
correspond to theta-hyperplanes of type 0 of $Y$, which do not exist by the previous argument.
Hence $X$ has no theta-hyperplanes of type 1.

Consider now a split curve $X$ of even genus and let $N$ be a node. Repeat the above construction: now
$Y$ has $2^{g-2}$ theta-hyperplanes of type 0, by the previous part of the proof,
hence $X$ has $2^{g-2}$ theta-hyperplanes of type 1 through $N$. Since $X$ has $g+1$ nodes,
we get
$$
t_1(X)=(g+1)2^{g-2}
$$
and we are done.

The rest of the proof is absolutely identical, and we leave it to the reader.
\qed

\

Finally, we have the following

\proclaim Lemma \mult.
Let $X\in V$; then a theta-hyperplane of type i appears with multiplicity $2^i$ as an irreducible
component of $\t (X)$.

\pf
This follows easily from [C], proof of (5.2), as explained in example (5.4).
\qed

\

\rk
The previous results imply that if $X$ and $Y$ are curves in $V$ such that $\t (X) = \t(Y)$, 
then they have the same number of nodes.

The reader can at this point check that
$\deg \t (X) = \sum 2^it_i(X)=2^{g-1}(2^g -1)$.

\

We conclude this section with a somewhat parenthetic discussion
 about cuspidal curves, from which the  rest of the paper is independent.
The goal is to describe the enumerative data of the configuration of theta-hyperplanes
for such curves. 

Let $Z\subset \P ^{g-1}$ be an integral curve 
of degree $2g-2$, having $\gamma$ ordinary cusps, that is,
$Z$ has $\g$ double points $q_1,....,q_\g$ such that if
$\nu :  \Zn \la Z$ is the normalization of $Z$, then $\nu ^{-1}(q_i) $ is one point
$p_i$. Notice that the genus of $\Zn$ is $h=g-\g$ and the arithmetic genus of $Z$ is
$g$. Let $K$ be the canonical class of $\Zn$; the projective model $Z$ is obtained
by a ($g-1$)-dimensional linear series in the  linear system
$|K + 2\sum _1^\g p_i|\cong \P ^{g+\g -2}$.

The discussion at the beginning of Section 2 can be extended to such cuspidal curves;
in particular we have:

\ni
(1) The definition of theta-hyperplanes, and that of  theta-hyperplanes of type $i$
($0\leq i\leq \g$) remains unchanged for such cuspidal curves.

\ni
(2) If $Z$ is generic, then $Z$ has a finite number of theta-hyperplanes.

\ni
(3) One can extend the definition of $\t (Z)$ to generic cuspidal curves, so that $\t
(Z)$ will be a non-reduced hypersurface of degree $N_g$, all of whose irreducible
components are theta-hyperplanes.

The details for $g=3$ 
are carried out in [CS].  Here we are only interested in the
enumerative aspects. 
  The theta-hyperplanes of $Z$ can be
stratified according to the number of cusps that they contain. Recall that $\t _i(Z)$
denotes the set of theta-hyperplanes containing exactly $i$ cusps, and that $t_i(Z) :=
\# \t _i(Z)$.

To better state and discuss the next result we need a new piece of notation.
A smooth curve of genus $h$ has $2^{2h}$ theta-characterstics. We shall denote by
$N^+_h:=2^{h-1}(2^h+1)$ the number of even theta-characteristics (recall that 
$N_h:=2^{h-1}(2^h-1)$ is the number of odd ones).

\proclaim Proposition \cusp. Let $Z\subset \P ^{g-1}$ be a general, integral curve of
degree
$2g-2$ having $\g$ ordinary cusps and no other singularities.
Then if $i<g-1$
$$
t_i(Z) = \cases{{\g \choose i}N_h & if \  $i \equiv \g $  (mod 2) \cr
{\g \choose i}N_h^+ & if \  $i \not\equiv \g $  (mod 2) \cr}
$$
 If $i = g-1$
$$ 
t_{g-1}(Z) =\cases{1, & if $\g = g-1$;\cr
g, & if $\g = g$.\cr}
$$

\

\pf
Let $S(Z^\nu)\subset \Pic ^{h -1} \Zn$ be the set of theta-characteristics of $\Zn$
(notation as above).
Denote by 
$\T (Z)$ the support of $\t (Z)$, in other words:
$\T (Z): =\cup _{i=0}^{\g} \t _i(Z)$
and, of course, $\# \T (Z) = \sum _1^\g t_i (Z)$.

In the first part of the proof we show that there is a surjective map of sets
$$
\tau : \T(Z) \la S(\Zn )
$$
all of whose fibers have cardinality $2^{\g -1}$. In particular, 
we get $\# \T (Z) = 2
^{g-2\g-1}$.

To define $\tau$, let $H\in \T(Z)$, then there exists a unique effective divisor
$D_H$ on $\Zn$ such that $\nu ^* (H\cap Z)=2D_H$. We have that $2D_H\sim K + 2\sum p_i$
and hence $D_H - \sum p_i$ is a theta-characteristic on $\Zn$.
We define
$$
\tau (H) = D_H - \sum _1^\g p_i.
$$
To show that $\tau$ is surjective, pick $T\in S(\Zn)$ and consider the divisor $D_T:=T+\sum
p_i$.
By genericity of $Z$ we have that $h^0(\Zn,D_T)=\g$;
there is a natural morphism
$$
\a : \P^{\g -1} \cong |D_T|\la |2D_T|=|K + 2 \sum p_i| \cong \P^{g+\g -2}
$$
defined by $\a (D)=2D$ (compare with the proof of Proposition \irr ).
Let now $W\subset |2D_T|$ be the $g-1$-dimensional linear subpace such that our
ambient space
$\P ^{g-1}$ is $W$.
The fiber $\tau ^{-1}(T)$ is given exactly by the points in $\Im \a \cap W$.

Now, $\a$ is the composition of the degree $2$ Veronese map with a general linear
projection. Thus $\dim \Im \a = \g -1$ and $\deg \Im \a = 2 ^{\g -1}$.
We obtain that $\Im \a \cap W$ is non-empty (by dimension count) hence $\tau$ is
surjective.
More precisely, our genericity assumption implies tha the intersection above is
transverse, that is, that $\Im \a \cap W$ is made of $2 ^{\g -1}$ distinct points.
Thus $\forall T \in S(\Zn)$ we have $\# \tau ^{-1}(T)= 2^{\g -1}$.

In the second part of the proof we work by induction on $\g$.
If $\g =1$ the description in the first part becomes particularly simple.

If $T$ is odd, that is, $T$ is effective
($\Zn$ is generic), the divisor $2(T + p_1)$ is cut out by a unique theta-hyperplane
of $Z$, passing through the cusp. 

If $T$ is even, that is, $T$ is not effective, we still have that $T+p_1$ is
effective (by Riemann-Roch) and that $|T+p_1|=\{ D_T\}$, where $D_T$ is an effective
divisor not containing $p_1$ in its support.
Hence $2D_T$ is cut out by a theta-hyperplane of type $0$ of $Z$.
Hence the case $\g =1$ of the Proposition is proved.

Let now $\g$ be arbitrary and $g-2>i>0$. 
Fix a set $G\subset \{q_1,....,q_\g \}$  made of $i$ cusps of $Z$ and let
$\pi : Z \la Z^G\subset \P ^{g-i-1}$ be the projection from the linear space spanned
by $G$, so that $Z^G$ is an integral curve of degree $2(g-i-1)$ having $\g -i$ cusps,
to which we can apply the induction hypothesis.
In particular
$$
t_0(Z^G) = 
\cases{N_h & iff \  $\g -i \equiv 0 $  (mod 2) $\Longleftrightarrow  \g \equiv i $
(mod 2)\cr  N_h^+ & iff \  $\g -i \not\equiv 0 $  (mod 2) $\Longleftrightarrow \g
\not\equiv i $ (mod 2)\cr}
$$
Moreover, $Z$ has precisely $t_0(Z^G)$ theta hyperplanes of type $i$ containing $G$;
since there are ${\gamma \choose i}$ choices for $G$ we get
$$
t_i(Z)={\gamma \choose i}t_0(Z^G)
$$
and we are done if $0<i<g-2$ .

Just like in the proof of Proposition \irr, the case $i=g-2$ needs to be treated
separately. Let $\pi : Z \la Z'\subset \Pl$  be the projection form a subset
of $G$
containing
$g-3$ cusps and let $Q'\in Z'$ be the image of the remaining cusp. $Z'$ is thus an
integral curve having $\g - g + 3$ ordinary cusps (one of which is $Q'$)
and no other singularities. The set of theta-hyperplanes of type $g-2$ of $Z$
containing $G$ is in one-to-one correspondence with the set 
$\t _{Q'}(Z')$ of theta-lines of type $1$ of $Z'$ containing $Q'$.
By induction
$$
\#\t_{Q'}(Z')=\cases{N_h & iff \  $\g -g+3 \equiv 1 $  (mod 2) $\Longleftrightarrow 
\g
\equiv g-2 $ (mod 2)\cr  
N_h^+ & iff \  $\g -g+3 \not\equiv 1 $  (mod 2)
$\Longleftrightarrow \g
\not\equiv g-2 $ (mod 2)\cr}
$$
Since of course $t_{g-2}(Z)={\g \choose g-2} \#\t_{Q'}(Z')$ we are done in this case.

The formula for $i=g-1$ is clear (the cusps are in general linear
position, by our genericity assumption). 

It remains to compute $t_0(Z)$.
Let us do it for $\gamma $ odd, for simplicity. It is obvious how to modify the
computation  that follows if $\g$ is even.
By the first part of the proof, every theta-characteristic of $\Zn$ determines $2^{\g
-1}$ theta -hyperplanes of $Z$. So far, we have found that the total number of
theta-hyperplanes corresponding to odd theta-characteristics of $\Zn$ is equal to
$$
t_\g (Z) + t _{\g -2}(Z) + ....+ t_3 (Z) + t_1 (Z) = N_h\cdot \sum _{i=1}^{\g +1
\over 2}{\g \choose 2i-1} =N_h \cdot2^{\g-1} 
$$
where the last equality follows from observing that $\sum _{i=1}^{\g +1
\over 2}{\g \choose 2i-1}$ is equal to the number of subsets having odd cardinality
of  a set of $\g $ elements; thus $\sum _{i=1}^{\g +1
\over 2}{\g \choose 2i-1} =  {2 ^{\g} \over 2}$. 

We conclude that
the odd theta-characteristics of $\Zn$ determine only theta-hyperplanes of type $i>0$.
In fact,
an analogous computation gives that the total number of theta-hyperplanes of type $i>0$
corresponding to even theta-characteristics of $\Zn$ is equal to 
$$
t_{\g-1} (Z) + t _{\g -3}(Z) + ....+ t_4 (Z) + t_2 (Z) = N_h^+\cdot \sum _{i=1}^{\g -1
\over 2}{\g \choose 2i} =N_h^+ \cdot(2^{\g-1} -1)
$$
since $\sum _{i=1}^{\g -1
\over 2}{\g \choose 2i}$ is the number of non-empty subsets of even cardinality in a
set with $\g$ elements.

By the first part of the proof we know that
$$
\# \T (Z) = 2^{\g -1} ( N_h^+ + N_h )
$$
and we just showed that
$$
\sum _{i=1}^\g  t_i(Z) = 2^{\g -1} (N_h^+ + N_h) - N_h^+.
$$
We conclude that $t_0(Z) =N_h^+$. \qed

\

\ni
{\bf 3. Recovering curves from their theta-hyperplanes.}
The first result here states that split curves are uniquely determined, among all curves in $V$,
by the weighted configuration of their theta-hyperplanes.
Using the terminology of [CS], we could restate the Theorem below by saying that split curves have the
theta-property.

\proclaim Theorem \splitthm. Let $X\in V$ be a  split curve, and let $Y\in V$ be
such that $\t (X) = \t (Y)$.
Then $X=Y$.

\pf
Let $X=C_1\cup C_2$ be a split curve and let $Y\in V$ be such that $\t (X)=\t (Y)$; denote $\T = \t (X)$.
Then $Y$ is also a split curve, since split curves can be characterized among all curves in $V$ as those
having either no hyperplanes of even type (if $g$ is even) or no hyperplanes of odd type (if $g$ is odd).
This property is clearly reflected on the structure of $\T$, whose irreducible components have all
multiplicity an odd power of $2$ if $g$ is even, or an even power of $2$ if $g$ is odd.

\

\ni
{\it Part 1: $X$ and $Y$ have the same nodes.} 
 Recall first that the nodes of $X$, being points on a rational
normal curve, are in general linear position, that is, any $h<g$ of them span a linear
space of dimension $h-1$. We shall now exhibit an algorithm to determine the nodes of
$X$ and
$Y$ as the $g+1$ ``cluster" points of the weighted configuration of hyperplanes given
by $\T$. More precisely, we shall prove that such nodes can be recovered from the
configuration of hyperplanes that they span. This is really a special case of a more
general phenomenon; see the discussion immediately following the proof 
of the Theorem.

For the  first step of the algorithm, pick a hyperplane $H\subset \T$ having
multiplicity
$2^{g-1}$, so that
$H$ contains exactly $g-1$ nodes of
$X$;  call such nodes $N_1,N_2,....,N_{g-1}$ (which we will now recover) and call the remaining nodes of
$X$  (not lying on $H$) $M_1, M_2$. 
There are ${g+1 \choose g-1} -1$ remaining theta-hyperplanes of type $g-1$; let $H'$ be one of them.
Of course 
$$
\dim H\cap H' = g-3
$$
and there are two possibilities:

\ni
(a) $H' = <N_{i_{1}},....,N_{i_{g-2}},M_j>$,
hence 
$$
J:=H\cap H' =<N_{i_{1}},....,N_{i_{g-2}}>.
$$
Observe also that such a $J$ can be obtained by intersecting $H$ 
with exactly 2 different $H'$, that is
$$
J= H\cap <N_{i_{1}},....,N_{i_{g-2}},M_1>=H\cap <N_{i_{1}},....,N_{i_{g-2}},M_2>.
$$
(b) $ H'=<N_{i_{1}},....,N_{i_{g-3}},M_1,M_2>$, 
then $K:=H\cap H'$ contains only $g-3$ nodes $N_i$'s and all such intersections are
distinct, since for every $g-3$-uple of nodes $N_i$s there is a unique $H'$ of type
(b) containing it.

We can distinguish intersections $H\cap H'$ of   type (a) (the $J$'s) from those of
type (b) (the $K$'s) using their multiplicity, since every $J$ is obtained from two
different
$H'$'s, while every $K$ is obtained in only one way.

For the second step, forget the $K$'s and, for a fixed $i$ with $1\leq i\leq
g-1$, denote
$J^i
\subset H$ 
$$
J^i= <N_1,....,N_{i-1},N_{i+1},....N_{g-1}>
$$
the $g-3$ dimensional subspace spanned by all $N_j$'s with the exception of $N_i$.
Of course $J^i$ is one of the intersections of type (a) above.

Now it is clear that the intersection of $J^i$ and $J^h$ is the $g-4$ linear space spanned by all $N_j$'s
with the exception of $N_i$ and $N_h$. Moreover, for any subset of distinct integers $
\{i_1,i_2,....,i_{g-2}\}\subset
\{1,2,....,g-1\}$ we get
$$
J^{i_1}\cap J^{i_2}\cap ....\cap J^{i_{g-2}} = N_h, \  \  \  \  \  \  h\not\in \{i_1,i_2,....,i_{g-2}\}.
$$
In this way we can recover all the nodes $N_1,....,N_{g-1}$ of $X$ and $Y$.

To determine the remaining two nodes of $X$ it is enough to repeat the argument starting from some
different
$H$. Or, go on with the proof of the Theorem, which only requires $X$ and $Y$ to have three nodes in
common.

\

\ni
{\it Part 2: Project from the nodes to show that $X=Y$.}
From now on, denote $N_1,....N_{g+1}$ the nodes of $X$ and $Y$, and let $Y= D_1\cup D_2$.

We use induction on $r$. The case $r=2$ is 4.1.1 in [CS].

Let $\pi _i$ be the projection from $N_i$ to $\P ^{r-1}$. Denote by $X_i$ and $Y_i$ the images
of $X$ and $Y$ via $\pi _i$. Thus $X_i$ and $Y_i$ are split curves in $\P ^{r-1} $ and
$\t(X_i) = \t (Y_i)$.
By induction $X_i = Y_i$.
We can choose two nodes (say $N_1$ and $N_2$)
such that 
$$\pi _1(C_1) = \pi _1 (D_1) = E_1$$
and
$$
\pi _2(C_1) = \pi _2 (D_1) = E_2$$
(and also, of course $\pi _1(C_2) = \pi _1 (D_2)$,  $\pi _2(C_2) = \pi _2 (D_2)$).
Let us show that $C_1=D_1$.
Denote by $S_i\subset \P ^r$ the cone over $E_i$ with vertex $N_i$
(that is 
$$
S_i = \cup _{P\in X}\overline{N_iP}= \cup _{Q\in Y}\overline{N_iQ} \  )
$$ 
then
$$
C_1\cup D_1 \subset S_1\cap S_2.
$$
We shall now see that this implies $C_1 = D_1$.
To do that, it is enough to show that the intersection
$ S_1\cap S_2$ is a scheme entirely supported on $C_1 \cap \overline{N_1N_2}$
(and thus it cannot contain a rational normal curve different from $C_1$).
By contradiction, assume that there is a point  $Q\in  S_1\cap S_2$ such that $Q$ does not belong to
$C_1$ or to the line joining $N_1 $ with $N_2$.
Then the line $\overline {N_iQ}$ must intersect $C_1$ in another point $P_i$. Thus we find 2 points,
$P_1$ and $P_2$ such that the plane $\Pi$ spanned by $N_1,N_2,Q$ satisfies
$$
P_1+P_2+N_1+N_2 \subset \Pi \cap C_1.
$$
Now this is not possible, in fact pick $r-3$ other points on $C_1$ and let $H$ be the hyperplane spanned
by $\Pi$ and those points; then 
$$
\deg H\cap C_1 \geq \deg \Pi \cap C_1 +r-3\geq r+1
$$ 
while 
$\deg C_1=r$ and
$C_1$ is not contained in a hyperplane.
(Notice that we did not need the points $P_1,P_2,N_1, N_2$ to be distinct)

We conclude that $C_1=D_1$.
In the same way we prove that $C_2 = D_2$ and hence $Y=X$, a contradiction.
\qed

\

\rk As suggested by the referee, what proved in Part 1 is a special case of the
following  more general result. Let $S=\{ P_1,....,P_t\}$ be a set of points in general
linear position in
$\P ^r$, with $t\geq r+1$. Let $\cal G$  be the configuration of 
${t \choose r}$ hyperplanes that they
span, that is: ${\cal G} = \{ H = < P_{i_1},...., P_{i_r}> \  {\rm with}\  \  1\leq
i_1<....<i_r\leq t\}$. Then $\cal G$ uniquely determines $S$.

Notice that the case $t=r+1$ is clear, and the next case $t=r+2$ is
precisely Part 1 of the proof of the Theorem. The general statement can be proved
using an analogous argument.
Let $H\in {\cal G}$ and rename the points in $S$ so that
$S = \{ N_1,....,N_r,M_1,....,M_s\}$ with $H= <  N_1,....,N_r>$ and $s+r=t$.
We shall show how to determine $N_1,....,N_r$. Let $H'\in {\cal G}$ and let us consider
$H\cap H'$. There are two posiibilities:

\ni
(a) $J :=H\cap H'$ is contained in $s+1$ elements of $\cal G$, that is, there exists
$H_1,....,H_s \in {\cal G}$ such that
$$
J= H\cap H_1=H\cap H_2=....=H\cap H_s
$$

\ni
(b) $H\cap H'$ is contained in less than $s+1$ elements of $\cal G$.

Now, (a) occurs for $H'$ containing at most one $M_i$ among its generators, that is,
for example 
$$J= H \cap <N_1,....,N_{r-1},M_1>=....= H \cap
<N_1,....,N_{r-1},M_s>=<N_1,....,N_{r-1}>$$

On the other hand, (b) occurs for $H'$ containing at least $2$ of the $M_i$.
In fact,
suppose that (say) $H' =<N_1,....,N_k,M_1,....,M_h>$ with $2\leq h \leq r$.
If $K=H\cap H' = H \cap H''$ then 
$H''=<N_1,....,N_k,M_{i_1},....,M_{i_h}>$, with $\{ M_1,....,M_h\}\cap
\{ M_{i_1},....,M_{i_h}\}= \emptyset$.
In fact if (say) $M_1 \in H''$, then $M_1 \in H'\cap H''$. Thus $M_1 \in H $ since 
$H'\cap H''\subset H\cap H' \cap H'$, and this is impossible. Thus
there exist at most ${s \over h}$ different hyperplanes $H'\in {\cal G}$
such that $K= H\cap H'$. Now we argue exactly as in the second step of the
algorithm used to prove  Part 1:  we recover $N_1,....,N_r$ by 
intersecting linear spaces  of codimension
$2$ of type (a), which are distinguished, being the only ones having multiplicity
$s+1$.

\

\ni
We conclude with  a result for  curves of genus $4$.

\proclaim Theorem \genus4. Fix $g=4$.
Let $X\in V$  be a general singular curve. If   $Y\in V$ is a curve such
that  there exists a quadric surface $Q$ containing $X$ and $Y$ and if $\t (X) = \t (Y)$,
then $X=Y$.

\rk
From a moduli theoretic point of view the condition that $X$ and $Y$ are contained in the same quadric does
not appear so restrictive, since all  quadrics of the same rank are projectively equivalent,
and every canonical curve in $\P ^3$ is contained in a  (unique) quadric.

\

\pf
First we show that $X$ and $Y$ have the same nodes. The fact that they have the same number of nodes is
an immediate consequence of Lemma \mult . Let $\T = \t (X)$.
Let us show that the nodes of $X$ are uniquely determined as the cluster points of $\t (X)$.

By Proposition \irr, if $X$ has a unique node $N$, then $\T$ contains $28$ planes 
$H_1,....,H_{28}$ of multiplicity $2$, passing through $N$ and tangent to $X$ at
smooth points (by  Proposition \irr ). The rest of $\T$ consists in $64$ distinct
planes not containing $N$. The node $N$ is determined as the unique point where all
the $H_i$'s intersect. In fact, if the intersection of the $H_i$  contained another
point  $M$, let $J$ be the line through $N$ and
$M$ and let $\phi : X^{\nu } \la \P ^1$ be the degree $4$ map  given by all planes through $J$
($ X^{\nu }$ is the normalization of $X$). Then the ramification locus of $\phi $ contains the divisors cut
on $X- N$ by the 28 planes $H_i$, contradicting the Riemann-Hurwitz formula ($ X^{\nu }$ has genus $3$).

If $X$ has 2 nodes, applying  Proposition \irr  we get that $\T$ contains 6 planes of multiplicity $4$,
which are the theta-planes of type $2$. Intersecting any two of them gives the line $J$ through the two
nodes. The intersection of $J$ with any plane of type $1$ (i.e., any plane appearing in $\T$ with
multiplicity $2$) determines one of the two nodes. Thus we recover the two nodes of $X$ which must
therefore coincide with the
two nodes of $Y$.

If $X$ has $3$ nodes, let $H$ be the plane containing all of them (the only component having multiplicity $8$
in $\T$). Intersecting $H$ with all the planes of multiplicity $4$,
each of which contains a pair of nodes
of $X$ (there are 12 such planes, 4 through each pair of nodes, by Proposition \irr ) 
gives 3 lines of $H$ forming a triangle of vertices the three nodes of $X$.

If $X$ has 4 nodes, then $\T$ contains $4$ planes of multiplicity $8$, each of which 
is spanned by the  three nodes
of $X$. The intersection of any three of them is thus a node of $X$ (and of $Y$).

If $X$ has more than 4 nodes, then $X$ is necessarily a split curve, which has already been treated
in Theorem \splitthm.

We now argue by contradiction,  assuming  that $X$ and $Y$ are distinct, irreducible 
(by Theorem \splitthm ) with $\d$ nodes.

Pick a node $N$ of $X$ and $Y$ and let $\pi : \P ^3 \la \P ^2$ be the projection from $N$.
Let $C=\pi (X)$ and $D =\pi (Y)$, so that $C$ and $D$ are irreducible plane quartic having 
the same $\d -1$ nodes.
We have
$$
\t (C) = \t(D)
$$
in fact the bitangent lines of $C$ and $D$ correspond to the theta-hyperplanes of type $1$
of $X$ and $Y$ containing $N$.
By [CS] (4.2.1, 5.1.1, 5.1.2, 5.2.1) we obtain that $C=D$.
Let now $S$ be the cone over $C$ of vertex $N$; by the above  analysis, $S$ contains both $X$ and $Y$.
But then
$$
X\cup Y \subset S\cap Q 
$$
which contradicts B\'ezout's Theorem, since $\deg X\cup Y = 12$, while $\deg S\cap Q= 8$.
\qed

\

\ni
We believe that the assumption that $X$ is singular can be removed from the statement 
of the Theorem,
and that a suitable analogue holds for curves of any genus.

\

\

\ni{\it Aknowledgments.} It is a pleasure to thank Edoardo Sernesi 
for introducing me to these
problems and for the  conversations we had about them. Special thanks to Emma
Previato for her enthusiastic encouragement and to the referee for important comments.
Part of this paper was written while I was visiting  the University of Strasbourg: 
my sincere gratitude to Olivier Debarre and
to the Department of
Mathematics for the invitation and for the wonderful
working conditions provided.

\

\

\ni
{\bf References }

\

\ni
[ACGH] E. Arbarello, M. Cornalba,  P. Griffiths, J. Harris: 
{\it Geometry of Algebraic Curves 1} 
Grundlehren der Mathematischen Wissenschaften, 267. Springer, New York-Berlin, 1985

\

\ni
[CS] L. Caporaso, E. Sernesi: {\it Recovering plane curves from their bitangents.}
Preprint Alg-Geom AG/0008239

\

\ni
[C] M. Cornalba: {\it Moduli of curves and theta-characteristics.}
Lectures on Riemann surfaces (Trieste, 1987), 560--589,
World Sci. Publishing, Teaneck, NJ, 1989

\

\ni {\bf Mathematics Subject Classification:} 14Nxx, 14Hxx.

\

\ni
L.Caporaso: Universit\`a degli Studi del Sannio, Benevento, Italy

and Massachusetts Institute of Technology, Cambridge, MA, USA

caporaso@math.mit.edu
\end